\documentclass[12pt]{article}
\usepackage{modif}
\usepackage{amsthm,
 }
\usepackage{graphicx}
\usepackage[colorlinks=true,citecolor=black,linkcolor=black,urlcolor=blue]{hyperref}
\newcommand{\doi}[1]{\href{http://dx.doi.org/#1}{\texttt{doi:#1}}}
\newcommand{\arxiv}[1]{\href{http://arxiv.org/abs/#1}{\texttt{arXiv:#1}}}

\theoremstyle{plain}
\newtheorem{theorem}{Theorem}
\newtheorem{lemma}{Lemma}[section]

\newtheorem{proposition}[theorem]{Proposition}

\theoremstyle{definition}

\theoremstyle{remark}

\def\3#1#2#3{_{#1=#2}^{#3}}   
\def\ca#1{\left|{#1}\right|}
\def\idn#1{\mathop\alpha(#1)}
\def\cqn#1{\mathop\omega(#1)}

\def\mval#1 {\delta\mathopen({#1}\mathclose)}
\def\mxval#1 {\Delta\mathopen({#1}\mathclose)}
\def\dist{\mathop{\rm dist}\nolimits}

\let\0=\emptyset

\def\lk{\mathop{{\rm lk}}\nolimits}

\def\e#1 {\ensuremath{\mathop{\rm
 e}\mathopen({#1}\mathclose)}}
\def\rr#1#2{\ensuremath{\mathop{\rm R}
 \mathopen({#1},{#2}\mathclose)}}

\def\br{\hfill\break}
\def\kk#1#2{\ensuremath{K_{{#1},{#2}}}}

\def\5{\ensuremath{{\cal W}_{13;1,5}}}

\def\v#1{\ensuremath{\deg(#1)}}
\def\vind#1#2{\ensuremath{\deg_{#1}(#2)}}
\def\vv#1#2{\ensuremath{\deg^{\,#1}(#2)}}
\def\vvind#1#2#3{\ensuremath{\deg_{#1}^{#2}(#3)}}

\def\kdots{\ifmmode,\ldots,\else,~$\ldots\,$, \fi}
\def\row#1#2{\ifmmode #1_1\kdots #1_{#2}\else
 $#1_1$\kdots$#1_{#2}$\fi}
\def\dash{\hbox{--}}

\long\def\clause#1#2{\par\smallbreak\hangafter=1\hangindent
 20pt\noindent{\hbox to20pt{{#1\hfill}}{#2}}\hfill\par
 \ifdim \lastskip <\smallskipamount \removelastskip
 \penalty 15\smallskip \fi}

\def\R{Ramsey }
\def\HoG{House of Graphs}

\title{\bf Edge number report 1: state of the art
estimates for $n \le 43$.}

\author{J\"orgen Backelin}

\begin{document}
 \maketitle

 \begin{abstract}
 This first extracted report contains all lower and
upper bounds for e-numbers \e3,k;n , for $n \le 43$,
that I know. All but 24 of them are known (exactly).
Very little of the proofs is given. A few
consequences for upper classical \R number bounds
are mentioned.
 %
%
\end{abstract}

 \section{Introduction.} \label{S:intr}
 Throughout the years, I have investigated
e-numbers, and updated my tables of these and of
properties for graphs with edge numbers close to the
respective e-number. The results have been collected
in the various updated versions of~\cite B. However,
that work is not easily accessible; not only since I
have not made it public, but since it is large, and
based on a somewhat complex terminology, both for
graph objects and for methods for dealing with them.

 At present, I'm integrating the consequences of
Goedgebeur's and Radziszowski's investigations
in~\cite{GR} into my tables. This is slow work; I
have now more or less finished it up to vertex
number 43. This has yielded a few improvements,
compared both to~\cite{GR} and to older versions
of~\cite B.

 I have received some criticism for not making my
results more accessible. In this report, I indeed
try to present the more recent ones, as regards
e-number bounds; but not the further Ramsey graph
properties. I believe that this makes it easier to
uaccess {\sl the conclusions\/}; but it makes it
harder to reproduce or improve {\sl the proofs}. I
outline a few proof examples; they may at least
illustrate the `Ramsey calculus' methods.

 Moreover I also discuss upper bounds
for e-numbers. This is an area not equally well
covered by the literature, I think, and I'm not sure
of how good the upper bounds I give here are,
compared to the state-of-the-art.

 Finally, the terminology is a bit experimentative.
I try to make it more conformant to other recent
state-of-the-art articles, and (against my
instincts) leave a good bit undefined. I'll be very
thankful for comments, both on this, and on the
factual content of this report.

\section{Definitions.} \label{S:def}
 Throughout this work, all graphs $G = (V,E)$ are
finite, simple, and undirected; and they are {\sl
triangle-free\/}; i.~e., the clique number $\cqn G
\le 2$.

 The {\sl second degree\/} of a vertex $v$ in a
graph $G$ is
 $$\vv2v = \vvind G2v := \sum_{w\in N(v)} \v w\,,$$
 where $N(v)$ is the set of vertices adjacent to
$v$. (The second degree is denoted $Z(v)$ in e.~g.\
\cite{GR}.) The induced $G$ subgraph on $V \setminus
(N(v) \cup \{v\})$ is denoted $G_v\,$.

 $G$ is an {\sl $(i,j;n,e)$-graph} and an {\sl
$(i,j;n)$-graph}, if $\cqn G < i$, its independence
number $\idn G < j$, $n(G) := \ca V = n$, and $e(G)
:= \ca E = e$.

For any positive integers $i$, $j$, and $n$, the
{\sl e-number\/} $\e i,j;n $ is the minimal number
$e$, such that there are $(i,j;n,e)$-graphs, or
$\infty$, if no $(i,j;n)$-graphs exist. They are of
great interest for finding improved bounds of \R
numbers
 $$\rr ij := \min \bigl(n : \e i,j;n = \infty
\bigr),$$
 but are also of interest in themselves.

 In this report, we only discuss the e-numbers
\e3,j;n . For the estimates, we shall use a few
linear or `piecewise linear' functions on two
integer variables, namely,

 $$f_1(n,k) = \max\, (0, n-k, 3n-5k, 5n-10k,
6n-13k)\,;$$
 $$f_2(n,k) = 8n-19.5k\,;$$
 $$f_3(n,k) = 9n-23k\,;\hbox{ and}$$
 $$f_4(n,k) = 6.8n-15.6k\,.$$

 Note, that $f_1(n,k) = 6n-13k$, if $n \ge 3k$.

 Occasionally, we mention the ``linear graph
invariant''
 $$t(G) := e(G)-6n(G)+13\idn G\,.$$

 \5 denotes the cyclic graph with 13 vertices
(conventionally named \row u{13}), and with two
vertices forming an edge if the absolute value of
their indices counted modulo 13 is either 1 or~5.
(This graph very often is denoted $H_{13}$.)

 For other concepts, background, et cetera, see the
bibliography. In particular, we shall discuss some
graphs given by means of {\sl extension patterns},
which provide recipies for constructing them
step-by-step; but neither the patterns and nor the
corresponding graphs are formally described here.

 \section{Known general values.}
 \label{S:kngen}

  For $n \le 3.25k+1.5$, all e-numbers are known.
(This indeed includes all \e3,k+1;n \ with $n \le
43$ and $k \ge 13$.) To begin with, we have
 \begin{proposition}
 \label{P:f1}
 For all positive integers $n$ and $k$,
 $$\e3,k+1;n \ge f_1(n,k)\,.$$
 The values are exact if and only if $n <
\rr3{k+1}$, and moreover either $n \le 3.25k-1$, or
$n=3.25n$.
 \end{proposition}

 For a proof, see e.~g.\ \cite{RK}. Note, that part
of the result is the fact that $t(G) \ge 0$ for all
(triangle-free) $G$.

 \begin{lemma}
 \label{L:s}
 Let $k$ and $n$ be positive integers, such that $3k
\le n < \rr3{k+1}$, but $\e3,k+1;n > f_1(n,k)$. Then
 $\e3,k+1;n = f_1(n,k)+1 \iff -1 < n-3.25k < 0$,
 $\e3,k+1;n = f_1(n,k)+2 \iff 0 < n-3.25k \le 0.5$,
 and $\e3,k+1;n \ge f_1(n,k)+3 \iff 0.5 < n-3.25k$,
 \end{lemma}

 The proof depends on deriving properties for graphs
with $t(G) \le 2$. In \cite{B}, indeed, all $G$ with
$t(G) \le 1$ are characterised, and sufficient
restrictions are found for those with $t(G)=2$.
(Actually, the complete characterising of the graphs
with $t(G)=0$ also is the main object of the
stand-alone manuscript \cite{B2}. The $t(G)=2$
result partly employs \cite{GR}.)

 Employing some constructions, we find that the
lower bound in the last part of lemma~\ref{L:s} is
exact in a few cases:

 \begin{lemma}
 \label{L:t3}
 If $3k \le n < \rr3{k+1}$ and $0.5 < n-3.25k \le
1.5$, then $\e3,k+1;n = f_1(n,k)+3$.
 \end{lemma}

 If $n > 3.25k + 1.5$, and moreover $k \le 12$, then
$\e3,k+1;n  > f_1(n,k)+3$; and I find it likely that
this should hold also for all higher $k$. Moreover,
I guess that
 \begin{equation}
 \label{E:conj}
 \e3,k+1;n \ge \max \bigl( f_2(n,k), f_3(n,k)
\bigr)\,,
 \end{equation}
 too; but I am far from being able to prove this.
The best general result I have for $n-3.25k\gg0$ is

 \begin{lemma}
 \label{L:c}
 For any $n$ and $k$,
 $$\e3,k+1;n  \ge f_4(n,k)\,.$$
 \end{lemma}
 (This is contained in \cite[proposition 13.5]B,
which is proved by means of a somewhat complicated
induction argument).

 \section{The other values for $n \le 34$.}
 \label{S:knsmall}

 For $n \le 34$, all \e3,k+1;n \ are known.
Actually, only 15 of them are `sporadic', i.~e., not
given by the known \R numbers, or in
section~\ref{S:kngen}; and they all have $n\ge22$
and $6 \le k \le 9$. Thus, they are included in the
following \e3,l;n \ table (where $l = k+1$):
 $$\begin{array}{|l|cccc|}
 \hline
 n\backslash l&7&8&9&10\\
 \hline
 22&60&42&30&21\\
 23&\infty&49&35&25\\
 24&\infty&56&40&30\\
 25&\infty&65&46&35\\
 26&\infty&73&52&40\\
 27&\infty&85&61&45\\
 28&\infty&\infty&68&51\\
 29&\infty&\infty&77&58\\
 30&\infty&\infty&86&66\\
 31&\infty&\infty&95&73\\
 32&\infty&\infty&104&81\\
 33&\infty&\infty&118&90\\
 34&\infty&\infty&129&99\\
 \hline
 \end{array}$$

Note, that all items under an $\infty$ in a column
also are $\infty$. In the sequel, in each column,
just the top $\infty$ (if any) is printed.

 \section{The other values and estimates for $35 \le
n \le 43$.}
 \label{S:estsmall}

 In the table, a single value indicates that this is
the exact e-value. Two values separated by a dash
($\dash$) are the best known lower and upper bounds
of the respective e-value. Again, $l = k+1$.

 $$\begin{array}{|l|ccccc|}
 \hline
 n\backslash l&9&10&11&12&13\\
 \hline
 35&140&107\dash108&84\dash85&68&55\\
 36&\infty&117\dash119&92\dash94&75&60\\
 37&&128\dash(132)&100\dash103&82&66\\
 38&&139\dash(143)&109\dash112&89\dash90&72\\
 39&&151\dash161&119\dash121&96\dash98&78\\
 40&&161\dash\infty&128\dash130&103\dash107&87\\
 41&&172\dash\infty&139\dash(150)&111\dash116&94\\
 42&&\infty&149\dash(160)&120\dash125&101\dash102\\
 43&&&159\dash(171)&129\dash134&108\dash111\\
 \hline
 \end{array}$$

 The upper bounds within parentheses are rather
preliminary; they are achieved by crude
constructions, made more or less on the fly, since I
am too ignorant to know where to look for the best
actually achieved upper bounds. I expect there to
have been constructions or computer enumerations
around for a while, giving better upper bounds for
all five or most of them.

 \section{Consequences for \R numbers.}
 \label{S:R}
 By hand calculations or by means of e.~g. the
matlab programme FRANK (\cite L)\footnote{The
version of FRANK that I employ includes a test for
raising the lower e-number bound in a few cases,
where the only formally possible degree
distributions all would have to contain either a
triangle of low-degree vertices, or a low-degree
vertex with too few low-degree neighbours (and thus
a too high second degree). In practice, this only
may happen, when the unraised e-number bound would
be close to, but slightly less than, the e-value for
some regular graph. This tweak yielded e.~g.\
$\e3,13;51 \ge 179$.}, it is fairly easy
to check for consequences for upper bounds on \R
numbers for any improvement of lower bounds of
e-numbers.
 As compared to the combined values from \cite{GR}
and older versions of \cite B, the sharper bounds
presented here yield just two improved upper \R
number bounds.

  It turned out that the improvement of the lower
bound for \e3,12;43 \ from~128 to~129 was crucial
for deducing that
 $$\rr3{19} \le 132,$$
 as reported in the latest dynamic survey on small
\R numbers (\cite{DS1}).

 The improvement of lower \e3,11;39 \ bound from~117
(\cite{GR}) to~119 suffices to prove that
 $$\rr3{16} \le 97\,.$$
 This bound is not (yet) included in the dynamic
survey.

 \section{A few proof hints.}
 \label{S:P}

 \subsection{Lower bounds.}
 Most of the `sporadic' lower bounds are found in
\cite{GR}; and/or are direct consequences of lower
bounds for smaller independence numbers. The
exceptions are the lower bounds for \e3,11;35 ,
\e3,12;38 , \e3,12;39 , \e3,13;41 , \e3,13;42 ,
\e3,12;43 , \e3,11;39 , and \e3,11;41 .

 The first six of these bounds, as well as the
`general' bounds, depend partly on theoretical
classification of some `lower' graphs, i.~e., graphs
with lower independence and vertex numbers;
likewise, the two last ones depend on computational
classification of some lower graphs. In all cases,
there is some use of properties deduced for some
lower graphs; and the general proof technique is to
assume the existence of a graph $G$ with
`offendingly' low $e(G)$, and then to deduce more
and more precise conditions for $G$, until finally a
contradiction is achieved. I'll provide a few
examples.

 \smallskip
 First, assume that $G$ is a (3,11;35)-graph with
$e(G) \le 83$; whence actually equality must hold.
We then successively may prove:

 \smallskip
 \indent$(a)\qquad \mval G > 2$;\br
 \indent$(b)\qquad \mval G > 3$;\br
 \indent$(c)$\qquad any vertex of degree 4 has at
most one neighbour of degree $\ge5$;\br
 \indent$(d)\qquad G_v$ has no \5 component for any
vertex of degree 5; and\br
 \indent$(e)$\qquad if $\v v=5$, then $\vv2v\le24$.

 \smallskip
 Property~$(a)$ is immediate from the \e3,10;n \
values.

 $(b)$ follows from~$(a)$, and from the fact that
any $(3,10;31)$-graph $H$ with $e(H) \le 74$ has
$\mval H \ge 2$, strictly if $e(H)=73$; and that
there are at most two vertices of degree 2 in $H$,
which (if indeed there are two of them) moreover
must be adjacent.

 $(c)$ is immediate from $(b)$, and the fact that
$\vv2v \le 17$ for any vertex of degree~4.

 $(e)$ is an immediate consequence of~$(d)$, and of
the fact that any $(3,10;29,58)$-graph does contain
a \5 component. On the other hand, $(e)$ directly
yields a contradiction, since it means that we could
calculate as if \e3,10;29 \ were at least~59.

 \smallskip
 This just leaves the deduction of~$(d)$ from~$(b)$
and~$(c)$, which is somewhat less immediate. Assume
for a contradiction that $\v v=5$, and that $G_v$
has a \5 component. Let $N(v) = \{\row w5\}$, and
let $U$ be the set of vertices in \5, which are not
adjacent to any $w_i$; in other words, $U = \{u \in
V(\5) : \vind Gu = 4\}$.

 Now, $\ca U \le 8$, since $U$ cannot contain an
independent 4-set; if it did, any edge between $U$
and $N(v)$ would be redundant (in the sense that
removing it from $G$ would leave a graph which also
did not contain an independent 11-set), but $G$ can
contain neither a redundant edge, nor a \5
component. Thus, and by inspection of \5, if $U$
were non-empty, then there were a $u_j \in U$ with
at most two neighbours in $U$, and therefore at
least two neighbours of degrees $\ge5$,
contradicting $(c)$.

 Thus, instead, $U = \0$; i.~e., each vertex in \5
is adjacent to at least one $w_i$. This makes it
possible to apply a ``decharging'' argument.
`Charge' each $u_j$ with a unit charge, 1; and then
`discharge' each $u_j$ by distributing its charge in
equal proportions to its $w_i$ neighbours. The total
charge after discharging must stay 13. However, no
$w_i$ can receive a charge larger than $2.5$; which
means that $N(v)$ in total cannot carry a higher
charge than $12.5$. This is a contradiction; which
indeed proves~$(d)$.

 \medskip
 For a second example, assume that $G$ is a
$(3,11;41)$-graph with $e(G)=138$. There are few
theoretic ways for such a graph  to be `realised
numerically'; in other words, if we let the degree
distribution (degree sequence) of the graph be
$(n_0, n_1\kdots n_{10})$, then there are just a
handful possible such sequences, for which the
resulting Graver-Yackel defect $\gamma(G)$ would be
non-negative (cf.\ \cite{GY} and~\cite{GR}). In
fact, also employing that a single vertex $v$ of
degree 8 would have $\vv2v \le 8\cdot7 = 56$, and
thus a positive defect, and repressing all leading
and trailing zeroes in the distributions, we would
have one of
 $$(11,30),\; (12,28,1),\; (1,9,31),\; (2,7,32),
\hbox{ and } (3,5,33)$$
 as degree distribution, with the total defect
$\gamma(G) = 3$, 1, 2, 1, and~0, respectively.

 Put $F := \{v \in V : \v v=7$ and $\vv2v = 48\}$.
In other words, $F$ is the set of non-defect
vertices of degree~7. Counting directly yields that
$\ca F \ge 27$, in each one of the cases.

 For any $f \in F$, $G_f$ is a $(3,10;33,90)$-graph.
Now, Goedgebeur and Radziszowski classified all
these graphs, and made a list of all 57099 of them
available on the {\sl\HoG\/} (\cite{GR}). Running
the NAUTY (\cite{MP}) command {~\tt countg~--Jd~} on
this list reveals that any such graph $H$ contains
an induced \kk24, and has $\mval H \ge 4$. Moreover,
a theoretical analysis shows that for any vertex $v$
with $5 \le \v v \le 7$, either $\mval G_v \ge 3$,
or $\mval G_v = 2$ and $\gamma(v)=3$, or $\gamma(v)
> 3$.

 Now, choose such an $f$; if there is a vertex $x$
of degree~8, actually choose $f \in F \cap N(x)$;
choose a $ \kk24 \subset V_f \subset V$, with
$V(\kk24) = \{a_1,a_2; \row b4\}$ and $\v{a_1} \le
\v{a_2} \le 7$, say. We now note, that
 $$\mval G_{a_i} \le \v a_{3-i} - 4, \hbox{ for }
i=1,2\,;$$
 and employ this in estimating the defects of the
$a_i$.

 If $\v{a_1}=5$, then $\gamma(a_2) \ge 4 > 3 \ge
\gamma(G)$, a contradiction. Likewise, if $\v{a_1} =
6$, then $\gamma(a_2) = 3$, whence then $\gamma(a_1)
= 0$; whence anyhow
 $$6 \le \v{a_1} \le 7 = \v{a_2}\,.$$
 If $\v{a_1}=7$, then both $a_1$ and $a_2$ are
defective, and the further defects in $G$ sum up to
at most~1, whence in particular then $\mxval G = 7$.
Moreover, if $\v{a_1}=7$, then not both $a_1$ and
$a_2$ may have defects $\ge2$, whence instead then
at least one of them has second valency 47, and thus
at least five neighbours of degree~7, of which at
least four belong to $F$. Thus, in this case, we may
assume that $f' := b_4 \in F$; while if $\v{a_1}=6$,
then let $f'$ be arbitrarily chosen in $F \cap \lk
a_2\,$. In either case, there is some \kk24 in
$V_{f'}$, and this would also carry a defect at
least~2, which would yield a total defect at least 4
in $G$, a contradiction.

 \subsection{Upper bounds.}
 For $n \le 4k = 4l-4n$ (but excepting $(n,l) \in
\{(17,6),\; (22,7),\; (27,8)\}$), there are
constructions, whose connected components either are
described by their extension patterns, or are one or
the other of two {\sl exceptional graphs\/}: The
cyclic graph \5 (the unique (3,5;13,26)-graph), and
the {\sl twisted tesseract\/} (a (3,6;16,32)-graph).
(The twisted tesseract also is denoted
$(2{\cal W}_{8;1,4})_{5i}$ in \cite B; i.~e., it
consists of two disjoint copies of ${\cal
W}_{8;1,4}$, with the $i$'th vertex in the first
copy connected to the $5i$'th one in the second copy
by an edge; where indices are taken modulo~8.)

 The extension pattern of a graph $G$ of the kind we
consider here includes a triangle free graph $T$,
such that
 $$e(T) \le 2n(T),$$
 $$\idn G = n(T),$$
 $$n(G) = 2n(T)+e(T), \hbox{ and}$$
 $$e(G) = n(T) + 2e(T) + {1\over2} \sum_{x\in V(T)}
\deg(x)^2\,.$$

 This yields that the graphs with only patterned
and/or exceptional graphs as components indeed
fulfil~(\ref{E:conj}). In fact, for `most' $k$ and
$n$ with $3.25k \le n \le 4k$, we have such graphs
realising equality in~(\ref{E:conj}). However, there
are some irregularities, for two reasons. First,
each \5 component contributes 4 to the independence
number of the graph; and there may not be an integer
number of such components that realises equality
in~(\ref{E:conj}). Second, in general, for a
connected patterned graph $G$ with $3.25 \idn G \le
n(G) \le 4\idn G$, equality only can be achieved by
having only vertices of degrees~3 and~4 in the
pattern graph $T$ (since other degree distributions
yield higher $\sum_{V(T)} \deg(x)^2$); which for
$(3,10;36)$-graphs would force the pattern graph to
be 4-regular, on~9 vertices. By inspection, there is
no such triangle-free graph; the closest possible
degree distribution is (2,5,2) vertices of degrees
(3,4,5), respectively.

 The upper bound 161 for \e3,10;39 \ is reported by
Goedgebeur and Radziszowski in \cite{GR}, where it
is noted that both they and Exoo have found huge
amounts of $(3,10;39,161)$-graphs $G$, but no
$(3,10;39)$-graph with a lower number of edges.

 For the five upper bounds within parentheses, let
$L$ be the regular $(3^8)$-type lace with constant
offsets (1,3), a $(3,9;32,104)$-graph. (Laces are
defined and investigated in \cite B; they form a
special class of patterned graphs.) Its family
$(\row v8)$ of apices consists of non-adjacent
vertices of degree 6, where moreover $\dist(v_i,v_j)
\ge 3$, if $i$ and $j$ have the same parity. The
upper \e3,10;37 \ (\e3,10;38 ) bounds are achieved
by a 4-extension (5-extension) of $H$, employing 3
(all~4) of the odd-indexed $v_i$, respectively; and
the upper \e3,11;41\hbox{---}43 \ bounds by making a
further extension of one of these, employing the
$v_i$ with even indices.

 \bibliographystyle{plain}

 \end{document}